# A Redesigned Benders Decomposition Approach for Large-Scale In-Transit Freight Consolidation Operations

Abdulkader S. Hanbazazah, Luis E. Abril, Nazrul I Shaikh[1], and Murat Erkoc

*Department of Industrial Engineering, University of Miami*

*Coral Gables FL USA*

February 2017

## Abstract

The growth in online shopping and third party logistics has caused a revival of interest in finding optimal solutions to the large scale in-transit freight consolidation problem. Given the shipment date, size, origin, destination, and due dates of multiple shipments distributed over space and time, the problem requires determining when to consolidate some of these shipments into one shipment at an intermediate consolidation point so as to minimize shipping costs while satisfying the due date constraints. In this paper, we develop a mixed-integer programming formulation for a multi-period freight consolidation problem that involves multiple products, suppliers, and potential consolidation points. Benders decomposition is then used to replace a large number of integer freight-consolidation variables by a small number of continuous variables that reduces the size of the problem without impacting optimality. Our results show that Benders decomposition provides a significant scale-up in the performance of the solver. We demonstrate our approach using a large-scale case with more than 27.5 million variables and 9.2 million constraints.

**Key Words:** *freight consolidation; third party logistics; mathematical programming; Benders decomposition.*

---

[1] Corresponding author. E-mail: n.shaikh@miami.edu

# 1. Introduction

Third party logistics (3PL) was a 157.2 billion dollar sector within U.S. that witnessed a growth rate of 7.4% in 2014 (3PL Market Analysis, 2015). The growth in online shopping is expected to further fuel this market and has caused a revival of interest in finding better solutions to the large scale in-transit freight consolidation problem, which allows the 3PL to cut costs by use economies of scale and reduction in package count (Xiaomin 2017). For a 3PL provider, the consolidation problem requires determining what shipments to consolidate at an intermediate gateway or terminal versus what to ship directly to a customer such that the overall shipment costs are minimized while the delivery time windows are honored. This paper and the solutions herein are motivated by our involvement with a major 3PL provider. We study their problem of in-transit consolidation of products being shipped from $n$ geographically disperse shippers to a single business customer via $m$ consolidation points. Each product has a pre-defined delivery time window that cannot be violated and each product is shipped to the intermediate consolidation point (called gateway) before being shipped to the destination. The decision is what gateway to send the products to and whether to consolidate the shipments into larger containers before shipping them to their end destination.

      Typically a business or corporate customer that employs 3PL have standing orders from multiple suppliers for multiple products across a planning horizon, where each product has a pre-specified shipment date and delivery time window. The 3PL providers pick up the products from the suppliers on given shipment dates and deliver to the customer within delivery time windows. All picked up products are first shipped to intermediate gateways before forwarded to the customer. A 3PL company usually has more than one gateway that provides flexibility pertaining to shipment costs and consolidation options. The routing decisions therefore need to be made for two legs: from suppliers to gateways and from gateways to the final customer. The first leg decision involves assigning the shipment to a particular gateway and selection of the transportation mode. The consolidation related decisions are made at the gateway. In the second leg, the carrier ships the products-either as consolidated shipments or as is so as to minimize shipment costs without violating the constraints set by the delivery time windows. When a shipment is not consolidated into a container, it is forwarded to the customer as individual shipment, which is typically more expensive.

We formulate the in-transit freight consolidation problem described above as a mixed-integer programming (MIP) problem and our key contribution is the development of a Benders decomposition based solution approach that provides a significant scale-up in the performance of the solver. The decomposition replaces a large number of integer "freight-consolidation" variables by a small number of continuous variables that reduces the size of the problem in terms of both the number of variables and constraint without impacting the optimality. Using our approach, we can solve to optimality a large-scale case with more than 27.5 million variables and 9.2 million constraints.

The remaining paper is organized into 5 sections. We present a brief literature review on the freight consolidation problem in Section 2. There are several useful research output reported in the literature on both the in-transit freight consolidation problem as well as Benders decomposition, and we point the reader to the relevant reviews. The proposed MIP model is presented in Section 3 and the Benders decomposition based reformulation of our model is presented in Section 4. Finally, a detailed case study elucidating the efficacy of the decomposition approach for solving large scale in-transit freight consolidation problem is presented in Section 5 with the conclusions and potential extensions discussed in Section 6.

## 2. Literature Review

The literature on in-transit freight consolidation problem is vast and includes various nuances that change the problem complexity as well as researchers focus. A thorough survey of this literature is presented by Guastaroba et al. (2016) who focus on the use of intermediate facilities in freight transportation planning and their application on three different settings: vehicle routing problems, transshipment problems, and service network design problems. One of the more comprehensive mathematical models for the in-transit freight consolidation problem was developed by Croxton et al. (2001). Their MIP formulation addresses some of the operational issues arising in merge-in-transit distribution systems. The model formulation accounts for various complex yet necessary features of an in-transit freight consolidation problem and includes the integration of inventory and transportation decisions, the dynamic and multimodal components of the application, and the non-convex piecewise linear structure of the cost functions. These two papers together give a good insight into the general setting of the problem as well as the specifics about the modeling and operationalization that the reader may look up for details.

Both papers establish that the in-transit freight consolidation problem is NP complete and hence researchers have been focusing on developing heuristic approaches to scale up as well as speed up the problem. Researchers have relied on (a) dual-based solution methods (Song et al. 2008), (b) column generation algorithms (Moccia et al. 2011; Dondo and Mendez 2014), (c) cutting-plane procedures and branch-and-bound heuristics (Croxton et al. 2003), (d) heuristic search (Popken 1994; Golias et al. 2012), (e) simulations (Qian and Xu 2012), and (f) decomposition based heuristics (Jin and Muriel 2009) to achieve the dual objective of scale-up and speed-up without compromising the quality of the solution. To the best of our knowledge, none of the papers in the domain have looked at Benders decomposition based approach to solve large scale in-transit freight consolidation problem. That apart, Fischetti et al. (2016) recently proposed a redesigned Benders decomposition for solving large scale MIP that uses a projected decision space for a "thinned out" version of the classic decision problem and show that the method enables significant scale-up and speed-up without impacting the optimality of the solution. The decomposition takes advantage of the new hardware and software technologies such as multi-core processors.

Our model builds on the model proposed in Croxton et al. (2001) along with the redesigned Benders decomposition approach proposed in Fischetti et al. (2016). We include two linear cost structures that correspond to shipment from shipper to the consolidation point and from the consolidation point to the customer respectively and a time constraint on each shipment in addition to the constraints accounted for by Croxton et al. (2001). Our setting is relevant to 3PL providers who need to solve the large scale in-transit freight consolidation problem on a frequent basis.

## 3. Model Formulation

In this section, we introduce the MIP model for the studied problem. The model tackles the case of in-transit consolidation of products being shipped from $n$ geographically spread shippers to a single final destination through $m$ gateways. Under a multi-period setting, products must be picked form the suppliers and routed to the destination within a given time window. Our model assumes that the freight from a supplier to the customer is divisible into different loads, which can be transshipped via different gateways and the 3PL provider may choose this option for consolidation opportunities. The customer has pre-determined pickup dates and due date windows for each product from suppliers that the 3PL provider is aware of and delivery deadlines are imposed as hard constraints.

The problem involves two stages. In the first stage, products are shipped from suppliers (shippers) located in different locations to one of several gateways such as ports. There are alternatives for the mode of transportation (usually land or air). For each transportation mode, cost of shipment is linear increasing in the amount of shipment. The transportation cost typically depends on the distance between the supplier and the gateway.

At gateways, the products are forwarded to the customers either as less-then-container-load (LCL) shipments or full-container-load (FCL) shipments. For the LCL option, the shipment cost is linear increasing in the shipment weight. If sufficient volume of products can be consolidated into a container without violating the delivery time-windows of the products, a more economic option of FCL shipment can be exercised. The goal of the models is to identify the optimal shipment routes and schedules over a planning horizon that minimizes the total transportation costs.

As mentioned earlier it is assumed that the shipments from the suppliers can be broken into pieces and routed to separate gateways on their way to the end delivery point, i.e., the customer. The carrier may choose the option of dividing the products picked up from a supplier into subsets, if doing so provides opportunities for FCL consolidation at the gateways. A shipment can be stalled at a gateway before it is moved to the second stage so that it can be coupled with other shipments and consolidated into a container. However, as mentioned earlier, products cannot be delayed beyond a certain point in time which results in late delivery. If they cannot be consolidated into a container in a timely fashion they must be forwarded as LCL shipments so as to make their respective delivery deadlines. Keeping products at the gateway incurs holding costs for the carrier, which is typically low in comparison to savings obtained from consolidation. FCL consolidation necessitates the introduction of integers variables that represent the number of containers used at each gateway in each period.

The model attempts to minimize the total cost over a set of multiple periods, $D$ (typically days in this context). It incorporates a set of shippers, $S$, a set of products, $P$, and a set of gateways, $H$. We denote $c1l_{s,h}$ and $c1a_{s,h}$ as the unit cost of shipment from supplier $s$ to consolidation gateway $h$ by land and by air respectively. Likewise, $c2_h$ is the cost of sending 1lb from gateway $h$ to the final customer, $c3_h$ is the cost of sending 1 container from gateway $h$ to the customer, and $ci_h$ is the inventory cost per lb realized by keeping the shipment at consolidation gateway $h$ for one time period. A time period in this context is typically a day. As such, in the rest of the paper we employ "day" as our time unit.

Decision variables for the model are as follows: $X_{p,s,h,d}$ is the weight in lbs of product $p$ sent from shipper $s$ to gateway $h$ on day $d$ by land, $Y_{p,s,h,d}$ is the weight in lbs of the items of product $p$ sent from shipper $s$ to gateway $h$ on day $d$ by air, $Z_{p,h,d}$ is the weight in lbs of the items of product $p$ sent from gateway $h$ to the final customer on day $d$ as a LCL shipment, $U_{p,h,d}$ is the weight in lbs of the items of product $p$ sent from gateway $h$ to the final customer on day $d$ as a FCL shipment, $T_{h,d}$ is the number of containers shipped from gateway $h$ on day $d$, $N_{p,d}$ is the total inventory in lbs of product $p$ delivered to the final customer on day $d$, and $I_{p,h,d}$ is the inventory of product $p$ in lbs at consolidation gateway $h$ on day $d$.

As for other parameters; $d_{p,s,h}$ is the weight of the items in lbs that must be picked up from shipper $s$ on day $d$, $k$ is the maximum capacity in lbs per container, $t1l_{s,h}$ is the number of days that a shipment takes by land from shipper $s$ to gateway $h$, $t1a_{s,h}$ is the number of days that a package takes by air from shipper $s$ to gateway $h$, $t2_{s,h}$ is the number of days it takes a package to ship from gateway $h$ to the final customer, $t_w$ is the length of the time window, and finally $D$ is the total number of periods (days) in the planning horizon.

The MIP model is then given below in Eqs. (1) – (8).

$$\min \sum_{\forall p,s,h,d} (c1l_{s,h} X_{p,s,h,d} + c1a_{s,h} Y_{p,s,h,d} + c2_h Z_{p,h,d} + c3_h T_{h,d} + ci_h I_{p,h,d}) \tag{1}$$

s.t.

$$\sum_{\forall h} X_{p,s,h,d} + \sum_{\forall h} Y_{p,s,h,d} = d_{p,s,d} \qquad \forall p, s, d \tag{2}$$

$$\sum_p U_{p,h,d} \leq k T_{h,d} \qquad \forall p, h, d \tag{3}$$

$$U_{p,h,d} + Z_{p,h,d} + I_{p,h,d+1} = \sum_{\forall s} X_{p,s,h,d-t1l_{sh}} + \sum_{\forall s} Y_{p,s,h,d-t1a_{s,h}} + I_{p,h,d} \quad \forall p, h, d \tag{4}$$

$$\sum_{\forall h} U_{p,h,d-t2_{s,h}} + \sum_{\forall h} Z_{p,h,d-t2_{s,h}} + N_{p,d} = \sum_{\forall s} d_{p,s,d-t_w} + N_{p,d+1} \quad \forall p, d: d \geq t_w \tag{5}$$

$$\sum_{\forall h: d \geq t2_{s,h}} U_{p,h,d-t2_{s,h}} + \sum_{\forall h: d \geq t2_{s,h}} Z_{p,h,d-t2_{s,h}} + N_{p,d} = N_{p,d+1} \quad \forall p, d: d < t_w \tag{6}$$

$$T_{h,d} \in \mathbb{N} \qquad \forall h, d \tag{7}$$

$$X_{p,s,h,d}, Y_{p,s,h,d}, Z_{p,h,d}, U_{p,h,d}, I_{p,h,d}, N_{p,d} \geq 0 \qquad \forall p, s, h, d \qquad (8)$$

The objective of the model is shown in Eq. (1), where we want to minimize the total shipping cost, which is composed of the fixed and variable costs of land freight and air freight from shippers to gateways, the cost of freight from gateways to the final customer broken into LCL and FCL shipments, and the cost of inventory held at the gateways. Eq. (2) ensures that scheduled pick-ups are carried out and shipped to gateways from a given supplier on a given day. Eq. (3) ensures that the amount of products shipped from gateway to the customer via containers does not exceed the capacity of the containers. Eq. (4) enforces that in-bound shipments, shipments on-hold and out-bound shipments are balanced at the gateways for a product type on a given day. Eq. (5) and Eq. (6) keep track of the flow balance at the customer site and ensure that the products are delivered to the customer by their due dates. Eq. (7) and Eq. (8) are the integrality and non-negativity constraints respectively.

The incorporation of the limits on the time windows enforces the feasibility constraints guaranteeing that the maximum time-span that a shipment may take from pick-up at any given supplier location to delivery at the final customer does not exceed certain time duration, $t_w$. As such, typically, consolidation of all products at gateways may not be possible across the time horizon. Moreover, the inclusion of the holding costs at the gateways may deter the storage of shipments until full truck load containers are completely loaded for shipment. As such, the optimal solution is typically a mix of individual LCL and FCL shipments. If the time windows are sufficiently large, air option is usually not utilized except for consolidation purposes at the gateways since they are usually much more expensive.

## 4. Benders Decomposition

Benders' decomposition is a method that is usually used for large mixed binary and integer optimization problems, where the problem is divided into smaller sub-problems which enable the global solution of the problem to be achieved. The model presented in Eq. (1)-(8) grows in size as $p, s, h, d$ increase and most of the real life 3PL in-transit consolidation problems cannot be solved to optimality on account of constraints on computational resources. We therefore tailor the Benders decomposition method with the "thin-out" approach presented by Fischetti et al.

(2016) to solve large scale in-transit freight consolidation problems. Our methodology involves two major parts.

In the first part, we obtain the linear programming (LP) relaxation of the original model by relaxing the integrality constraints given in Eq. (7). The solution of the LP relaxation model provides us with a lower bound for the problem, which we can use as a starting point for the overall implementation. Next, we use the principles of Benders decomposition to find an upper bound solution for our problem. We adapt its solution philosophy to our problem under the following considerations:

a. Let us organize our objective function in following two parts: $FO = f(X) + g(Z) + h(T)$, where $f(X)$ is the cost for the first leg of the transshipments, that is cost of shipping products from suppliers to gateways. Here, $X$ is the array of all decision variables of the first leg. Let $g(Z)$ involve the costs related to the LCL shipments from gateways to the customer and holding inventory at the gateways. Finally, $h(T)$ is the part of the objective function that captures the FCL shipment costs from gateways to the customer.

b. Under Benders' strategy, when we fix T, our integer variable, the problem left to solve is of LP class. Under this view, we can rewrite our problem:

$$\text{Min } q(T) + h(T) \tag{9}$$
$$\text{st.} \qquad T \in \mathbb{N} \tag{10}$$

Here $q(T)$ is the solution to the following problem:

$$\text{Min } f(X) + g(Z) \tag{11}$$

st.

$$\sum_p U_{p,h,d} \leq k\bar{T}_{h,d} \qquad \forall p, h, d \tag{12}$$

$$\sum_{\forall h} X_{p,s,h,d} + \sum_{\forall h} Y_{p,s,h,d} = d_{p,s,d} \qquad \forall p, s, d \tag{13}$$

$$U_{p,h,d} + Z_{p,h,d} + I_{p,h,d+1} = \sum_{\forall s} X_{p,s,h,d-t1l_{sh}} + \sum_{\forall s} Y_{p,s,h,d-t1a_{s,h}} + I_{p,h,d} \qquad \forall p, h, d \tag{14}$$

$$\sum_{\forall h} U_{p,h,d-t2_{s,h}} + \sum_{\forall h} Z_{p,h,d-t2_{s,h}} + N_{p,d} = \sum_{\forall s} d_{p,s,d-t_w} + N_{p,d+1} \qquad \forall p, d: d \geq t_w \tag{15}$$

$$\sum_{\forall h: d \geq t2_{s,h}} U_{p,h,d-t2_{s,h}} + \sum_{\forall h: d \geq t2_{s,h}} Z_{p,h,d-t2_{s,h}} + N_{p,d} = N_{p,d+1} \qquad \forall p, d: d < t_w \tag{16}$$

$$X_{p,s,h,d}, Y_{p,s,h,d}, Z_{p,h,d}, U_{p,h,d}, I_{p,h,d}, N_{p,d} \geq 0 \qquad \forall p, s, h, d \tag{17}$$

In the above model, $\bar{T}$ is a given integer value rather than a decision variable. Note that if the above model is unbounded for some $\bar{T} \in \mathbb{N}$, then the mathematical model given in Eqs. (9) – (10) is also unbounded, which in turn implies unboundedness of the original problem. If the model defined by Eqs. (11) – (17) is bounded, than we can obtain $q(T)$ by solving its dual. Furthermore, assuming feasibility of the region of the dual, we can enumerate all extreme points $(\alpha_p^1, \dots, \alpha_p^I)$ and extreme rays $(\alpha_r^1, \dots, \alpha_r^J)$. Notice that by solving our problem for $q(T)$, we can also access its dual $\alpha_p^i, \alpha_r^j$ variables. This implies that the mathematical model defined by Eq. (11) – (17) can be viewed as a sub problem. Let $q$ represent the optimal objective function value of this subproblem. Consequently, our master problem becomes:

$$Min\ q + h(T) \tag{18}$$

$$\text{st.} \qquad \left(\alpha_p^i\right)'(b - BT) \leq q \qquad \forall i = 1, \dots, I \tag{19}$$

$$\left(\alpha_r^j\right)'(b - BT) \leq 0 \qquad \forall i = 1, \dots, J \tag{20}$$

$$T \in \mathbb{N} \tag{21}$$

Here Eq. (18) represents Benders' optimality cut, and Eq. (19) represents Benders' feasibility cut, where B is a parameter matrix whose elements come from the coefficients of all constraints that involve the integer variable $T$ and b is a vector whose elements are the parameters from the coefficients of constraints in Eqs. (13) – (16).

Given that there exist an exponential number of extreme points and extreme rays of the dual of $q(T)$, generating all constraints of the type of Eq. (19) and Eq. (20) is not practical. Instead, we solve our Benders decomposition starting with a subset of these constraints and solving a relaxed master problem which yields a candidate solution. We iterate solving the sub problem and the master problem until the bounds meet, i.e., $q$ converges to a value.

## 5. Case Study

In this section, we introduce a case study that illustrates the implementation our suggested solution methodology discussed in the previous section. The case is based on a real life problem that a major 3PL provider in the United States faces frequently. In the case study considered here, the customer is a manufacturer of generic drugs based out of Puerto Rico and it provided the 3PL company with the supply data-i.e., product details, quantity, shipping date, shipping location, delivery time window for one calendar year. There are a total of 722 products originating from 104 supply locations spread over 25 states in the USA. The descriptive statistics of these products and a summary of the shipping locations are presented in Tables 1 and 2 respectively. The expected delivery pattern, i.e., the quantity to be picked up from the supply location on a specific day is presented in Figure 1. Products are aggregated based on their weights rather than volume since the latter one is relatively insignificant.

**Table 1.** Descriptive statistics of supply data

|  | Min | 1st Qu. | Median | Mean | 3rd Qu. | Max. | Obs. |
|---|---|---|---|---|---|---|---|
| Products (lbs.) | 15 | 417.2 | 1414 | 3028 | 4520 | 40000 | 722 |
| Daily Shipment quantity (lbs.) | 0 | 0 | 2054 | 5975 | 9051 | 76537 | 365 |

**Table 2.** Number of scheduled pick-ups across states and 365-day time horizon

| Origin State | Total | Origin State | Total | Origin State | Total |
|---|---|---|---|---|---|
| AL | 1 | MA | 27 | PA | 63 |
| AZ | 19 | MD | 3 | SD | 9 |
| CA | 98 | MN | 3 | TN | 1 |
| DL | 5 | MO | 10 | TX | 22 |
| FL | 8 | NC | 24 | UT | 5 |
| GA | 48 | NJ | 37 | VA | 1 |
| IL | 91 | NM | 5 | WI | 59 |
| IN | 72 | NY | 33 | Grand Total | 722 |
| KY | 37 | OH | 41 | | |

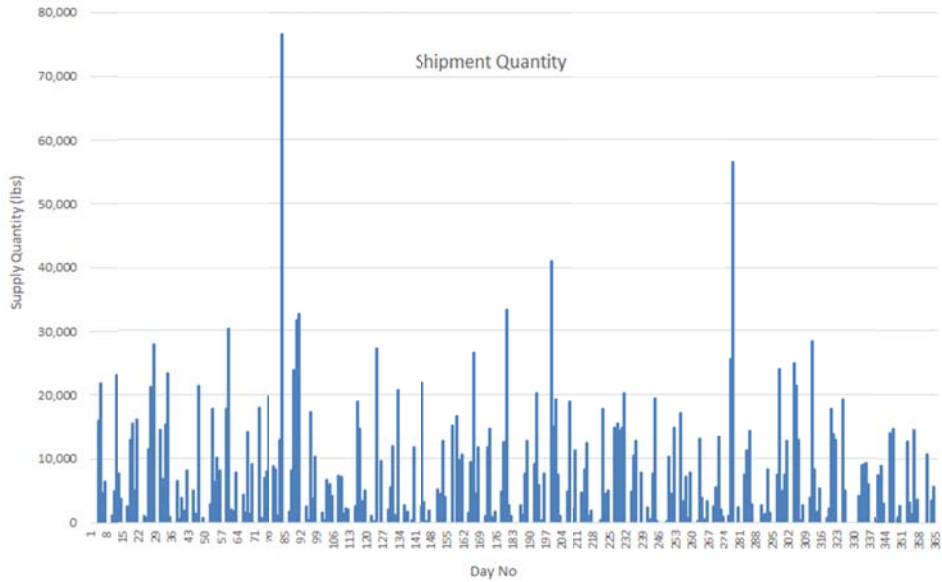

**Figure 1.** Demand for shipment across the 365-day time horizon

The customer requires pick-up dates from the suppliers and strict constraints on tardiness; all the products need to reach the manufacturer's site set in Puerto Rico within a 9-day window. As pointed out above, the 3PL company's objective is to determine the lowest cost at which it could satisfy the customer's shipment requirements. The 3PL company operates three gateways on the East Coast of the USA which are located at Port of Elizabeth, NJ; Port of Jacksonville, FL; and Port of Miami, FL. In the first leg, products are collected from the supplier locations summarized in Table 2 and shipped to one of these gateway locations using trucking. At gateways, the shipments are consolidated into LCL and FCL shipments and forwarded to San Juan, PR via ocean freight. The local delivery of the products in San Juan is omitted from the problem since the impact of this stage to the overall problem is negligible.

The unit transportation costs in the first stage depend on the location of the supplier and where the product will be shipped to. The costing at the 3PL company is done based on zones. Locations in USA are allocated to these zones for costing purposes. The zones are shown in Figure 2. The zone matrix given in Table 3 is used to identify the transportation times and unit transportation costs to gateways given in Table 4. We note that Zones 9 and 12 are used as destination zones since the gateways are located in these zones. The zoning structure and cost values are slightly modified to protect 3PL company's private information.

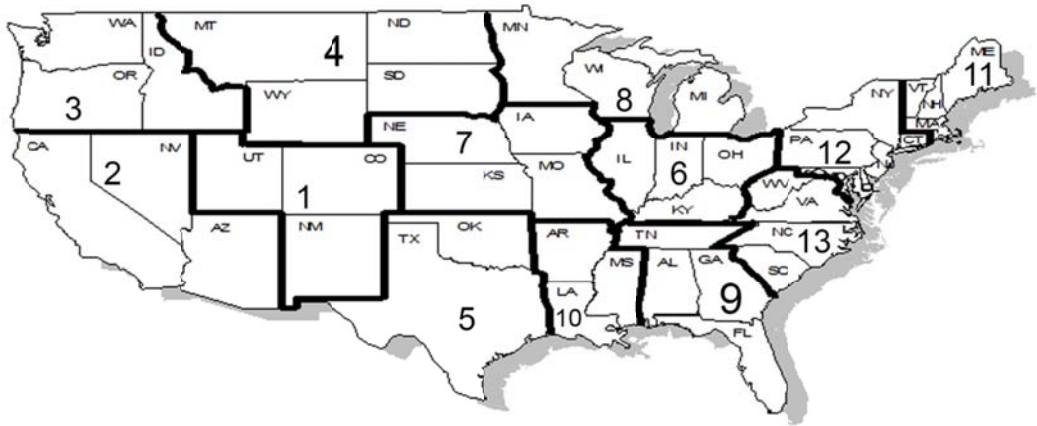

**Figure 2.** Zoning within the US and Zone Matrix

**Table 3.** The zone matrix

| FROM ZONE | TO ZONE 1 | 2 | 3 | 4 | 5 | 6 | 7 | 8 | 9 | 10 | 11 | 12 | 13 |
|---|---|---|---|---|---|---|---|---|---|---|---|---|---|
| 1 | A1 | B2 | C1 | C1 | B3 | C2 | B2 | C2 | D3 | C3 | D3 | D3 | D4 |
| 2 | B2 | A1 | B1 | D2 | B3 | C2 | C2 | C2 | D4 | D4 | D4 | D4 | D4 |
| 3 | C1 | B1 | A1 | C1 | D2 | D2 | D2 | D2 | D4 | D4 | D4 | D4 | D4 |
| 4 | C1 | D3 | C1 | B2 | E1 | D1 | D3 | C2 | E2 | E2 | E2 | E2 | E2 |
| 5 | B3 | B3 | D2 | E1 | A1 | B2 | B2 | B3 | B3 | B2 | C3 | C3 | C3 |
| 6 | C2 | C2 | D2 | D1 | B2 | A1 | B1 | B1 | B2 | B3 | B3 | B3 | B3 |
| 7 | B2 | C2 | D2 | E | B2 | B1 | A1 | B1 | C2 | B3 | C3 | C3 | C2 |
| 8 | C2 | C2 | D2 | C2 | B3 | B1 | B1 | A1 | B3 | C3 | B3 | B3 | B3 |
| 9 | D3 | D4 | D4 | E2 | B3 | B2 | C2 | B3 | A1 | B1 | B3 | B2 | B1 |
| 10 | C3 | D4 | D4 | E2 | B2 | B3 | B3 | C3 | B1 | A1 | C3 | C3 | B2 |
| 11 | D3 | D4 | D4 | E2 | C3 | B3 | C3 | B3 | B3 | C3 | A1 | A2 | B2 |
| 12 | D3 | D4 | D4 | E2 | C3 | B3 | C3 | B3 | B2 | C3 | A2 | A1 | B2 |
| 13 | D4 | D4 | D4 | E2 | C3 | B3 | C2 | B3 | B1 | B2 | B2 | B2 | A1 |

Table 4. Shipment times and unit costs between zones

| | A1 | A2 | B1 | B2 | B3 | C1 | C2 | C3 | D1 | D2 | D3 | D4 | E1 | E2 |
|---|---|---|---|---|---|---|---|---|---|---|---|---|---|---|
| **Shipment Time (In Days)** | | | | | | | | | | | | | | |
| Time | 2 | 2 | 3 | 3 | 3 | 4 | 4 | 4 | 5 | 5 | 5 | 5 | 6 | 6 |
| **Variable Cost Per Pound (In Cents)** | | | | | | | | | | | | | | |
| Cost per Pound | 29 | 33 | 33 | 37 | 41 | 37 | 41 | 44 | 37 | 41 | 44 | 46 | 44 | 46 |

The company's goal is to gain cost advantage by consolidating multiple products into containers at gateway locations before shipping them to Puerto Rico. The cost and time-length of shipments from the three gateways to Puerto Rico are presented in Table 5.

Table 5. LCL and FCL Costs

| Ports | Time (Days) | LCL Rate Per 100 lbs. | Container Capacity (lbs.) | FCL Rate Per Container | Threshold |
|---|---|---|---|---|---|
| Jacksonville, FL | 1 | $25.50 | 48,000 | $4,773.00 | 39.0% |
| Elizabeth, NJ | 2 | $17.13 | 48,000 | $4,805.46 | 58.4% |
| Miami, FL | 1 | $16.02 | 48,000 | $3,888.00 | 50.6% |

The 3PL provider's problem requires making decisions on what products to ship through what consolidation center so that the transportation cost from supplier to the manufacturing unit is minimized while all the due date constraints are met. Clearly, this problem can be modeled using the MIP introduced in Section 3, where there are 722 items, 104 supply points, 3 gateways, and 365 time periods. We note that typically, there is also a fixed cost for pick up at the supplier site in the first stage. However, since this cost is fixed and identical across all locations ($80 per pick-up in the case study) and it applies before the shipments are split to gateways, they do not affect the optimality of the solution obtained by the proposed MIP model. The resultant MIP model has more than 27.5 million variables and 9.2 million constraints. Attempts were made to solve the MIP using CPLEX on an Intel (R) Xeon (R) CPU Es-268 WO @ 3.10 GHz (dual processor) with 64 GB RAM machine; however, all attempts at solving the MIP failed on account of lack of sufficient computational resources. Subsequently, the solution methodology presented in Section 4 was applied and results are discussed next.

**Results**

In order to solve the in-transit merging optimization problem using Benders decomposition approach, we divided the problem into integer and linear parts. The linear part consists of the delivery of the packages sent from the shippers to the gateways, while the integer parts consisted of the merger of products at the consolidation stations and their shipment using FCL containers, as well as individual shipments (using LCL containers). The total cost is then the result of adding the individual values of the three cost components: the cost of freight from suppliers to the gateways (the linear part of the model), the cost of freight from the gateways to the clients using FCL containers, and the cost of freight from the gateways to the clients using LCL containers. The fixed costs of the pick-ups at the first stage are added to the solution of the model so as to find the overall annual cost. The results obtained for this case are summarized in Table 6.

Solving the linear relaxation of the model yielded the following results: The linear part of the model had a total cost of $657,399.67; the FCL part of the model resulted in a total cost of $163,469.83; and the LCL part had a cost of $0. It is straightforward to see that the relaxed problem allocates all shipments of the second stage to FCL containers since the unit cost is lower and fractional numbers for containers are allowed due to LP relaxation. This resulted in a total cost of $876,629.50.

**Table 6.** LP relaxation and Benders decomposition results

| Problem | No. Containers | Shipping Costs | | | |
|---|---|---|---|---|---|
| | | f(x) | | g(t,z) | Total |
| | | fix cost | variable cost | | |
| LP Relaxation | 0 | $55,760.00 | $657,399.67 | $163,469.83 | $876,629.50 |
| Benders decomposition (delivery exactly on the 9th day) | 6 | $55,760.00 | $657,400.00 | $235,000.82 | $948,160.82 |
| Benders decomposition (delivery within 9 days) | 13 | $55,760.00 | $660,207.00 | $207,534.15 | $923,501.15 |

When implementing Benders decomposition to solve the in-transit merging problem, we first consider the scenario, where the customer expects a delivery exactly 9 days after a pick-up. Occasionally, early delivery is regarded as inconvenience by the client since they schedule the pick-up dates based on just-in-time production and avoid carrying input inventory. We capture this case by simply removing the variable $N_{p,d}$ from the proposed MIP model. In this case, we obtain an optimal objective value of $948,160.82. Overall, only 6 FCL containers were possible under this scenario. At the end, we observed that about 75% of the costs were incurred in the first stage in this case.

When we allow early delivery (i.e., replace the constraint of delivering on a specific day by a more relaxed constraint of delivering within a specific time window), more consolidation alternatives become feasible and this leads to lower costs. In our case, applying the Benders decomposition approach to the problem that allows early deliver yields a total cost of $923,501.15 which represents a $25,000 reduction in total costs on account of 7 additional FCL (13 in total) consolidations.

An additional advantage from the schedule provided by our model is that it ensures that all of the deliveries are carried out within the time window. This is a significant improvement for the company who delivered about 20% of the shipments outside the delivery time window. Their shipment time performance is depicted in Figure 3.

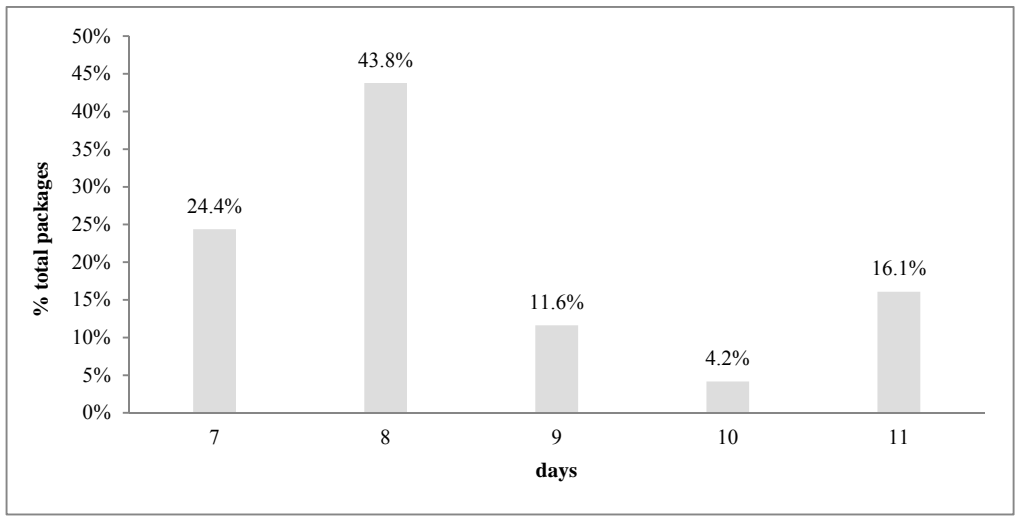

**Figure 3.** 3PL Actual Data - Distribution of delivery days for all products

The proposed optimization model which allows early delivery produced a solution with delivery performance depicted in Figure 4. The solution suggests a more uniform distribution in terms of delivery times. Approximately 25% of the shipments are consolidated into FCL shipments at gateways in the suggested solution. We believe that providing customer satisfaction by guaranteeing timely deliveries is paramount in the freight industry and the implementation of our proposed models ensure that highest quality service can be provided by the 3PL company.

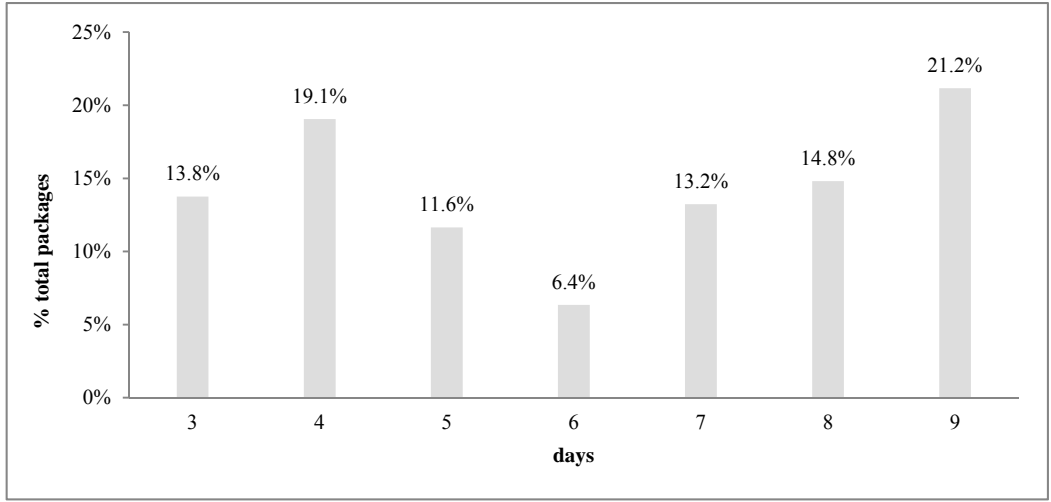

**Figure 4.** Distribution of delivery days for all products

## 6. Conclusions and Future Work

In this paper, we developed a MIP model that considers the case of in-transit consolidation of products being shipped from multiple shippers to a single business customer via multiple gateways that serve as consolidation points. The shipments have pre-specified pick up dates with delivery time windows across a multi-period time horizon. The problem is composed of two legs. In the first leg, products are shipped from suppliers to gateway locations, where shipment cost is a linear function of the package weight and distance between the supply point and the selected gateway. The shipments are forwarded from gateways to the customer's site either using LCL or FCL. The latter one is the cheaper option with lower unit costs however it is possible only if sufficient amount of shipments from the first stage can be consolidated at a given gateway. The delivery time windows impose constraints on consolidation opportunities since products must be delivered before their respective deadlines.

Due to complexity of the problem, the proposed model cannot be used to solve realistic size instances in its monolithic form. To facilitate practical use of the model, we propose a decomposition approach adapted from Benders decomposition method where the large numbers of integer "freight-consolidation" variables are replaced by a small number of continuous so as to reduce the size of the problem without impacting the optimality. Using a case study adopted from real life application, we showed that Benders decomposition provides a significant scale-up in the performance of the solver and we can solve a large-scale case with more than 27.5 million variables and 9.2 million constraints to optimality. Thus, the proposed redesigned Benders decomposition based approach solves large-scale in-transit freight consolidation problems optimally and efficiently.

Our solution has several practical benefits as well. The implementation of such a method will not only reduce the total costs for the 3PL providers, but will also enable them to solve larger-problems. In future work, we plan to extend the scope of our model to multiple customers that potentially facilitates more consolidation options at gateways.

# Appendix

Sets:

| | Description |
|---|---|
| $D$ | Set of time periods |
| $S$ | Set of shippers |
| $H$ | Set of gateways |
| $P$ | Set of products |

Model parameters:

| | Description |
|---|---|
| $c1l_{s,h}$ | Cost of sending 1 lbs from supplier $s$ to gateway $h$ by land |
| $c1a_{s,h}$ | Cost of sending 1 lbs from shipper $s$ to gateway $h$ by air |
| $c2_h$ | Cost of sending 1 lbs from gateway $h$ to the final customer |
| $c3_h$ | Cost of sending 1 container from gateway h to the final customer |
| $ci_h$ | Inventory cost per lbs at gateway $h$ per period |
| $d_{p,s,h}$ | Weight of product $p$ in lbs sent from shipper $s$ on day $d$ |
| $k$ | Maximum capacity in lbs per container |
| $t1l_{s,h}$ | Number of days that a package takes by land to arrive from shipper $s$ to gateway $h$ |
| $t1a_{s,h}$ | Number of days that a package takes by air to arrive from shipper $s$ to gateway $h$ |
| $t2_{s,h}$ | Number of days it takes a package to arrive from gateway $h$ to the final customer |
| $t_w$ | Length of the time window |

Decision Variables:

| | Description |
|---|---|
| $X_{p,s,h,d}$ | Weight in lbs of product $p$ sent from shipper $s$ to gateway $h$ on day $d$ by land |
| $Y_{p,s,h,d}$ | Weight in lbs of product $p$ sent from shipper $s$ to gateway $h$ on day $d$ by air |
| $Z_{p,h,d}$ | Weight in lbs of product $p$ sent from gateway $h$ to the final customer on day $d$ |
| $I_{p,h,d}$ | Weight of the of product $p$ in lbs at gateway $h$ on day $d$. |
| $U_{h,d}$ | Weight in lbs sent from gateway $h$ to the customer at day $d$ using a container |
| $N_{p,d}$ | Weight of the inventory (items delivered early) in lbs at the final customer on day $d$ |